\documentstyle[times,aps,amsthm,amsmath,amsfonts,%
    enumerate,twocolumn,floats]{revtex}

\input epsf

\def\nullset{\phi}

\def\We{We}
\def\we{we}
\def\our{our}

\newcommand{\EKLPa}{1}
\newcommand{\EKLPb}{2}
\newcommand{\fulmek}{3}
\newcommand{\MRR}{4}
\newcommand{\Macdo}{5}
\newcommand{\Sagan}{6}
\newcommand{\Stana}{7}
\newcommand{\Stanb}{8}
\newcommand{\Turn}{9}
\newcommand{\crummycite}[1]{[#1]}

\newcommand{\symmdiff}{{\bigtriangleup}}
\newcommand{\Prob}{{\rm Prob}}

\begin{document}

\wideabs{
\title{Domino tilings with barriers \\ \ \\ {\it In memory of Gian-Carlo Rota}}
\author{James Propp$^1$ and Richard Stanley$^2$}
\address{$^1$University of Wisconsin, Madison, WI 53706}
\address{$^2$Massachusetts Institute of Technology, Cambridge, MA 02139}
\maketitle

\begin{abstract} In this paper, \we~continue the study of domino-tilings of
Aztec diamonds (introduced in \crummycite{\EKLPa} and \crummycite{\EKLPb}). In particular, \we~look
at certain ways of placing ``barriers'' in the Aztec diamond, with the
constraint that no domino may cross a barrier. Remarkably, the number of
constrained tilings is independent of the placement of the barriers. \We~do not
know of a simple combinatorial explanation of this fact; \our~proof uses the
Jacobi-Trudi identity.
\end{abstract}}

\medskip

{\sf (NOTE: 
This article has been published in the Journal of Combinatorial Theory,
Series A, the only definitive repository of the content that has been
certified and accepted after peer review.  Copyright and all rights
therein are retained by Academic Press.  
You can also access this article via IDEAL (the International Digital
Electronic Access Library) at 
\begin{center}
{\tt http://www.idealibrary.com} 
\end{center}
or
\begin{center}
{\tt http://www.europe.idealibrary.com}\ .
\end{center}
This material may not be copied
or reposted without explicit permission.)
}

\vspace{0.2in}

{\bf I. Statement of result.}

\vspace{0.2in}

An {\bf \boldmath Aztec diamond of order $n$}
is a region composed of $2n(n+1)$ unit squares,
arranged in bilaterally symmetric fashion
as a stack of $2n$ rows of squares,
the rows having lengths
$$2, 4, 6, ..., 2n-2, 2n, 2n, 2n-2, ..., 6, 4, 2.$$
A {\bf domino} is a 1-by-2 (or 2-by-1) rectangle.
It was shown in [\EKLPa] that
the Aztec diamond of order $n$ can be tiled by dominoes
in exactly $2^{n(n+1)/2}$ ways.

Here we study {\bf barriers},
indicated by darkened edges of the square grid
associated with an Aztec diamond.
These are edges that no domino is permitted to cross.
(If one prefers to think of a domino tiling of a region
as a perfect matching of a dual graph
whose vertices correspond to grid-squares
and whose edges correspond to pairs of grid-squares
having a shared edge,
then putting down a barrier in the tiling
is tantamount to removing an edge from the dual graph.)
Figure 1(a) shows an Aztec diamond of order 8 with barriers,
and Figure 1(b) shows a domino-tiling
that is compatible with this placement of barriers.

\begin{figure*}
\begin{center}
\leavevmode
\epsfbox[0 0 360 156]{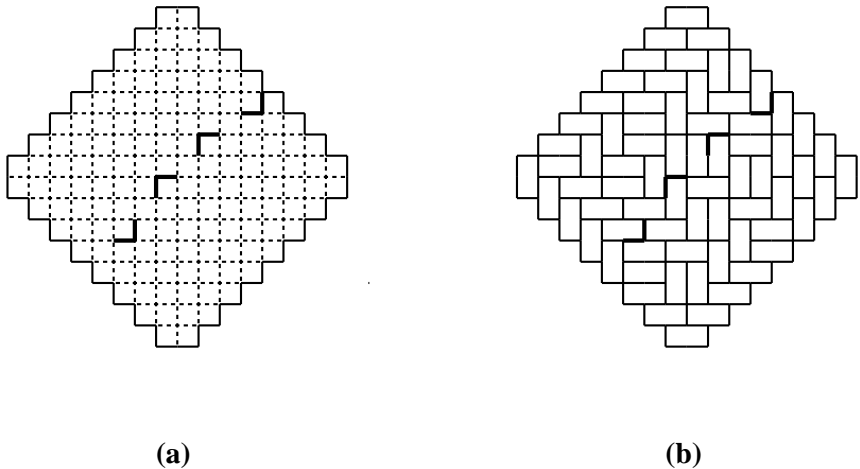}
\end{center}
\caption{Barriers and tiling.}
\label{fig-bar}
\end{figure*}


%

The barrier-configuration of Figure 1(a) has special structure.
Imagine a line of slope 1 running through the center of the Aztec diamond
(the ``spine''), passing through $2k$ grid-squares, 
with $k=\lceil n/2 \rceil$.
Number these squares from lower left (or ``southwest'')
to upper right (or ``northeast'')
as squares 1 through $2k$.
For each such square, 
we may place barriers on its bottom and right edges (a ``zig''),
barriers on its left and top edges (a ``zag''),
or no barriers at all (``zip'').
Thus Figure 1 corresponds to the sequence of decisions
``zip, zig, zip, zag, zip, zag, zip, zig.''
Notice that in this example, for all $i$,
the $i$th square has a zig or a zag if $i$ is even
and zip if $i$ is odd.
Henceforth (and in particular in the statement of the following Theorem)
we assume that the placement of the barriers has this special form.

\vspace{0.1in}
\noindent
{\bf Theorem 1:}
Given a placement of barriers in the Aztec diamond as described above,
the number of domino-tilings compatible with this placement is
$2^{n(n+1)/2}/2^k$.
\vspace{0.1in}

Some remarks on the Theorem:

(1) The formula for the number of tilings makes no mention at all
of the particular pattern of zigs and zags manifested by the barriers.
Since there are $k$ even-indexed squares along the spine, there are
$2^k$ different barrier-configurations, all of which are claimed to have
equal numbers of compatible tilings.

(2) Each domino-tiling of the Aztec diamond is compatible with exactly one
barrier configuration (this will be explained more fully in section II).
Hence, summing the formula in the Theorem over all barriers configurations,
one gets $2^k \cdot 2^{n(n+1)/2} / 2^k$, which is $2^{n(n+1)/2}$, the total
number of tilings.

(3) 180-degree rotation of the Aztec diamond switches the odd-indexed and
even-indexed squares along the spine, so the Theorem remains true if we consider
barrier-configurations in which the $i$th square has a zig or a zag if $i$
is odd and zip if $i$ is even.

\vspace{0.2in}

\noindent
{\bf II. Preliminaries for the proof.}

\vspace{0.2in}

Consider a particular tiling of an Aztec diamond,
and consider a particular square along the spine.
If that square shares a domino with the square to its left, or above it,
then placing a zag at that square is incompatible with the tiling.
On the other hand, if the square shares a domino with the square to its right,
or below it, then placing a zig at that square is incompatible with the tiling.
It follows that for each domino-tiling, there is a unique compatible way of
placing zigs and zags along the spine.  This holds true whether one only
puts zigs and zags at every other location along the spine (as in Figure 1(a))
or at every location along the spine.  In the case of the tiling depicted in
Figure 1(b), the full sequence of zigs and zags goes ``zag, zig, zig, zag, zig,
zag, zag, zig.''

Each such sequence must contain equal numbers of zigs and zags.
For, suppose we color the unit squares underlying the Aztec diamond
in checkerboard fashion, so that the squares along the spine are white
and so that each white square has four only neighbors (and vice versa).
The barriers divide the Aztec diamond into two parts, each of which
must have equal numbers of black and white squares (since each part can
be tiled by dominoes).  It follows that the white squares along the spine
must be shared equally by the part northwest of the diagonal and the
part southeast of the diagonal.

Given a sequence of $k$ zigs and $k$ zags, let
$1 \leq a_1 < a_2 < ... < a_k \leq 2k$ be the sequence of locations of the zigs,
and let
$1 \leq b_1 < b_2 < ... < b_k \leq 2k$ be the sequence of locations of the zags.
Note that the sets $A = \{a_1,...,a_k\}$ and $B = \{b_1,...,b_k\}$
are disjoint with union $\{1,2,...,2k\}$.
Let us call them a {\bf balanced (ordered) partition} of $\{1,2,...,2k\}$.
It is proved in section 4 of [\EKLPa]
that the number of compatible domino-tilings of the Aztec diamond of order $n$
is
\begin{equation}
\left( \prod_{1 \leq i < j \leq k} \frac{a_j - a_i} {j-i} \right)
\left( \prod_{1 \leq i < j \leq k} \frac{b_j - b_i} {j-i} \right)
2^{k'(k'+1)}
\end{equation}
where $k' = \lfloor n/2 \rfloor$.
(This is equivalent to Theorem 2 in [\MRR].)
For instance, the tiling shown in Figure 1(b) determines the balanced partition
$A = \{2,3,5,8\}$, $B = \{1,4,6,7\}$, and there are 2025 compatible tilings.

If we sum (1) over all balanced partitions of $\{1,2,...,2k\}$
we must of course get $2^{n(n+1)/2}$.  Theorem 1 claims that if we sum (1) over
only those balanced partitions $A,B$ which have certain specified even numbers 
in $A$ (and the remaining even numbers in $B$), we get $2^{n(n+1)/2}/2^k$.
Thus, to prove Theorem 1, it suffices to prove
that
\begin{equation}
\sum_{(A,B)}
\left( \prod_{1 \leq i < j \leq k} \frac{a_j - a_i} {j-i} \right)
\left( \prod_{1 \leq i < j \leq k} \frac{b_j - b_i} {j-i} \right)
\end{equation}
is independent of $A^* \subseteq \{2,4,...,2k\}$,
where the $(A,B)$ in the sum ranges over all balanced partitions of 
$\{1,2,...,2k\}$
such that $A \cap \{2,4,...,2k\} = A^*$.
Note that in this formulation, $n$ has disappeared from the statement of
the result, as has the Aztec diamond itself.

\vspace{0.2in}

\noindent
{\bf III. Restatement in terms of determinants.}

\vspace{0.2in}

We can interpret the left-hand side of (2) using Schur functions
and apply the Jacobi-Trudi identity.
The expression
\begin{equation}
\prod_{1 \leq i < j \leq k} \frac{a_j-a_i}{j-i}
\end{equation}
is equal to the number of semistandard Young tableaux of shape
$\lambda=(a_k-k, a_{k-1}-(k-1), ..., a_2-2, a_1-1)$ with parts at most $k$.
That is to say, if one forms an array of unit squares
forming left-justified rows of lengths $a_k-k,...,a_1-1$
(from top to bottom),
(3) gives the numbers of ways of filling in the boxes with numbers
between $1$ and $k$ 
so that entries are weakly increasing from left to right
and strictly increasing from top to bottom.

For background information on Young tableaux, Schur functions, and
the Jacobi-Trudi identity, see [\Macdo], [\Sagan], [\Stana], or [\Stanb].
In particular, for the definition of Schur functions and a statement
of the Jacobi-Trudi identity,
see formulas (5.13) and (3.4) of [\Macdo],
Definition 4.4.1 and Theorem 4.5.1 of [\Sagan],
or 
Definition 7.5.1 and Theorem 7.11.1 of [\Stanb].

If we associate with each semistandard Young tableau the monomial
$${x_1}^{m_1} {x_2}^{m_2} \cdots {x_k}^{m_k}$$
where $m_i$ is the number of entries equal to $i$ in the tableau,
then the sum of the monomials associated with the tableau
is the Schur function $s_\lambda(x_1,x_2,...,x_k,0,$
$0,...)$.
By the Jacobi-Trudi identity,
this is equal to the determinant
\[
\left| \begin{array}{cccc}
        \mbox{$h_{a_k-k}$} & \cdots & \mbox{$h_{a_k-2}$} & \mbox{$h_{a_k-1}$} \\
        \mbox{$h_{a_{k-1}-k}$} & \cdots & \mbox{$h_{a_{k-1}-2}$} & \mbox{$h_{a_{k-1}-1}$} \\
        \vdots & & \vdots & \vdots \\
        \mbox{$h_{a_1-k}$} & \cdots & \mbox{$h_{a_1-2}$} & \mbox{$h_{a_1-1}$}
        \end{array} \right|
\]
where $h_m$ is 0 if $m<0$ and otherwise is equal to 
the sum of all monomials in $x_1,...,x_k$ with total degree $m$
(so that $h_0 = 1$, $h_1 = x_1+...+x_k$,
$h_2 = x_1^2+x_1x_2...+x_k^2$, etc.).

Thus, if we let $v(m)$ denote the length-$k$ row-vector
\[
(h_{m-k}, \cdots, h_{m-2}, h_{m-1}) ,
\]
we see that the summand in (2) is the determinantal product
\[
\left| \begin{array}{c} v(a_k) \\ \vdots \\ v(a_1) \end{array} \right|
\cdot
\left| \begin{array}{c} v(b_k) \\ \vdots \\ v(b_1) \end{array} \right|
\]
specialized to $x_1 = x_2 = ... = x_k = 1$.
To prove Theorem 1, it will suffice to show that this product,
summed over all balanced partitions $(A,B)$ with 
$A \cap \{2,4,...,2k\} = A^*$,
yields
\[
\left| \begin{array}{c} v(2k) \\ v(2k-2) \\ \vdots \\ v(2) \end{array} \right|
\cdot
\left| \begin{array}{c} v(2k-1) \\ v(2k-3) \\ \vdots \\ v(1) \end{array} \right|
.
\]
For, since this expression is independent of $A^*$,
and since the sum of this expression
over all $2^k$ possible values of $A^* \subseteq \{2,...,2k\}$
is $2^{n(n+1)/2}$
(by the result proved in [\EKLPa]),
the value of the expression must be $2^{n(n+1)/2-k}$,
as claimed in Theorem 1.

It is interesting to note that
one can also evaluate the preceding determinantal product directly.
Appealing to the Jacobi-Trudi identity,
we see that the product is
$$s_\sigma (x_1,...,x_k) s_\tau (x_1,...,x_k)$$
where $\sigma = (k,k-1,...,1)$
and $\tau = (k-1,k-2,...,0)$.
It is known that
\[ s_\sigma (x_1,...,x_k) = x_1 \cdots x_k \prod_{i<j} (x_i+x_j) \] 
and
\[ s_\tau (x_1,...,x_k) = \prod_{i<j} (x_i+x_j) \ ,\] 
so that the determinantal product is
\[ x_1 \cdots x_k \prod_{i<j} (x_i+x_j)^2 . \] 
Setting $x_1=...=x_k=1$, we get $2^{k(k-1)}$.
Multiplying this by the factor $2^{k'(k'+1)}$ from (1),
we get $2^{k(k-1)+k'(k'+1)}$.
It is simple to check that regardless of whether $n$ is even or odd,
the exponent $k(k-1)+k'(k'+1)$ is equal to $n(n+1)/2-k$,
as was to be shown.

\vspace{0.2in}

\noindent
{\bf IV. Completion of proof.}

\vspace{0.2in}

We can deduce the desired identity as a special case of a general
formula on products of determinants.  This formula appears as
formula II on page 45 (chapter 3, section 9) of [\Turn],
where it is attributed to Sylvester.  
However, we give our own proof below.

Suppose $(A^*,B^*)$ is a fixed partition of $\{2,4,...,2k\}$ into
two sets, and let $v(1),...,v(2k)$ be {\it any} $2k$ row-vectors of length $k$.
Given $A \subseteq \{1,...,2k\}$ with $|A| = k$,
let $||A||$ denote the determinant of the $k$-by-$k$ matrix
\[
\left| \begin{array}{c} v(a_1) \\ v(a_2) \\ \vdots \\ v(a_k) \end{array} \right|
,
\]
where $A = \{a_1,a_2,...,a_k\}$
with $a_1 < a_2 < ... < a_k$.
Abusing terminology somewhat, we will sometimes think of $A$
as a set of vectors $v(a_i)$, rather than as a set of integers $a_i$.

\vspace{0.1in}
\noindent
{\bf Theorem 2:}
\begin{equation}
\sum_{(A,B)} ||A|| \cdot ||B|| = ||\{1,3,...,2k-1\}|| \cdot ||\{2,4,...,2k\}||
\label{e1} \end{equation}
where $(A,B)$ ranges over all balanced partitions of $\{1,2,...,2k\}$ with 
$A \cap \{2,4,...,2k\}$ $=A^*$,
$B \cap \{2,4,...,2k\}$ $=B^*$.
\vspace{0.1in}

Remark: This yields as a corollary the desired formula of the last
section, with an extra sign-factor everywhere to take account of the
fact that we are stacking row-vectors the other way.

Proof of Theorem 2: Every term on the left is linear in
$v(1),...,v(2k)$, as is the term on the right; hence it suffices to
check the identity when all the $v(i)$'s are basis vectors for the
$k$-dimensional row-space.

First, suppose that the list $v(1),...,v(2k)$ does {\it not} contain
each basis vector exactly twice.  Then it is easy to see that every
term vanishes.

Now suppose that the list $v(1),...,v(2k)$ contains each basis vector
exactly twice.  There are then $2^k$ ways to partition $\{1,...,2k\}$
into two sets $A,B$ of size $k$ such that $||A|| \ ||B|| \neq 0$,
since for each of the $k$ basis vectors we get to choose which copy
goes into $A$ and which goes into $B$.  However, not all of these
partitions occur in the sum on the left, since we are limited to
partitions $(A,B)$ with $A \supseteq A^*$, $B \supseteq B^*$.  Call
such balanced partitions {\bf good}.

Suppose that the basis vectors $v(1),v(3),...,v(2k-1)$ are not all
distinct; say $v(s)=v(t)$ with $s,t$ odd, $s<t$.  Then, for every good
balanced partition $(A,B)$ that makes a non-zero contribution to the
left-hand side, we must have $s \in A$ and $t \in B$ or vice versa.
But then $(A \symmdiff \{s,t\}, B \symmdiff \{s,t\})$ (where
$\symmdiff$ denotes symmetric difference) is another good balanced
partition.  We claim that it cancels the contribution of $(A,B)$.
For, if one simply switches the row-vectors $v(s)$ and $v(t)$, one
introduces $t-s-1$ inversions, relative to the prescribed ordering of
the rows in the determinant; specifically, each $i$ with $s<i<t$ has
the property that $v(i)$ is out of order relative to whichever of
$v(s),v(t)$ is on the same size of the new partition.  Ordering the
row-vectors properly introduces a sign of $(-1)^{t-s-1} = -1$.  This
leads to cancellation.

Finally, suppose that $v(1)$, $v(3)$, $...$, $v(2k-1)$ are all
distinct, as are $v(2)$, $v(4)$, $...$, $v(2k)$.  We must check that
the sole surviving term on the left has the same sign as the term on
the right.  This is clear in the case where $A^* = \{1,3,...,2k-1\}$
and $B^*$ is empty, for then the two terms are identical.  We will
prove the general case by showing that if one holds $v(1),...,v(2k)$
fixed while varying $(A^*,B^*)$, the sign of the left side of the
equation is unaffected.  For that purpose it suffices to consider the
operation of swapping a single element from $A^*$ to $B^*$.  Say this
element is $s$ (with $s$ odd).  Then there is a unique $t \neq s$
(with $t$ even) such that $v(s)=v(t)$.  Let us swap $s$ with $t$ in
the term on the left side of the equation; since $v(s)=v(t)$, the
determinants are not affected.  In performing the swap, we have
introduced either $t-s-1$ (if $t>s$) or $s-t-1$ (if $s>t$) inversions,
relative to the prescribed ordering of the rows.  Since both
quantities are even, we may re-order the rows in the determinants so
that indices increase from top to bottom, without changing the sign of
the product of the two determinants.  We now recognize the modified
term as the sole non-vanishing term associated with $(A^* \setminus
\{s\}, B^* \cup \{s\})$.  Since this term has the same sign as the
term associated with $(A^*,B^*)$, and since the sign is correct in the
base case $(\{1,...,2k-1\},\nullset)$, the correctness of the sign for
all partitions of $\{1,...,2k-1\}$ follows by induction.

This concludes the proof of Theorem 2, which in turn implies
Theorem~1.~$\ \Box$ 

\textsc{Remark.} An identity equivalent to summing both sides of
equation (\ref{e1}) for all $2^k$ sets $A^\ast$ (with row $v(m)$
specialized to $(h_{m-k},\dots,h_{m-2},h_{m-1})$ as needed for
Theorem~1) is a special case of an identity proved combinatorially by
M. Fulmek [\fulmek] using a nonintersecting lattice path argument. It
is easily seen that Fulmek's proof applies equally well to prove our
Theorem~2. Thus Fulmek's paper contains an implicit bijective proof of
Theorem~2.

\vspace{0.2in}

\noindent
{\bf V. Probabilistic application.}

\vspace{0.2in}

One can define a probability distribution on ordered partitions
of $\{1,$ $2,$ $...,$ $2k\}$ into two sets of size $k$,
where the probability of the partition $(A,B)$ is
\[
\frac{
\left( \prod_{1 \leq i < j \leq k} \frac{a_j - a_i} {j-i} \right)
\left( \prod_{1 \leq i < j \leq k} \frac{b_j - b_i} {j-i} \right)
}{2^{k^2}} .
\]
Theorem 1 is equivalent to the assertion that the $k$ random events
$2 \in A$, $4 \in A$, $...$, $2k \in A$ are jointly independent,
and it is in this connection that it was first noticed.
As a weakening of this assertion, we may say that the events
$s \in A$ and $t \in A$ are uncorrelated with one another
when $s$ and $t$ are both even (or both odd, by symmetry).

\vspace{0.1in}
\noindent
{\bf Theorem 3:}
For $1 \leq m \leq k$, let $N_m$ be the random variable
$|A \cap \{1,...,m\}|$, where $(A,B)$ is a random partition of $\{1,...,2k\}$
in the sense defined above.
Then $N_m$ has mean $m/2$ and standard deviation at most $\sqrt{m/2}$.
\vspace{0.1in}

Proof: Define indicator random variables
\[
I_i = \left\{ \begin{array}{ll}
        1 & \mbox{if $i \in A$,} \\
        0 & \mbox{if $i \in B$,} \\
        \end{array} \right.
\]
so that $N_m = I_1 + I_2 + ... + I_m$.  Each $I_i$ has expected value $1/2$,
by symmetry, so the expected value of $N_m$ is $m/2$.  To estimate the
variance, split up the terms of $N_m$ into $N_m^{\rm odd} = I_1 + I_3 + ...$
and $N_m^{\rm even} = I_2 + I_4 + ...$.  The terms in each sum are independent
random variables of variance $\frac{1}{4}$, so the variance of $N_m^{\rm odd}$ 
is $\frac{1}{4}\lceil m/2 \rceil$ and the variance of $N_m^{\rm even}$ 
is $\frac{1}{4}\lfloor m/2 \rfloor$.  It follows from the Cauchy-Schwarz
inequality that the standard deviation of $N_m$ is at most
$\sqrt{\frac{1}{4} \lceil m/2 \rceil}
+\sqrt{\frac{1}{4} \lfloor m/2 \rfloor}
\leq \sqrt{m/2}$,
as was to be shown.

\vspace{0.1in}

The significance of the random variables $N_m$ is that 
(up to an affine renormalization)
they are values of the 
``height-function'' associated with a random domino-tiling
of the Aztec diamond (see [\EKLPa]).
Theorem 3 tells us that if one looks along the spine,
the sequence of differences between heights of consecutive vertices
satisfies a weak law of large numbers.

\vspace{0.2in}

\noindent
{\bf VI. Open problems.}

\vspace{0.2in}

One open problem is to find a combinatorial (preferably bijective)
proof of Theorem 1.
For instance, one might be able to find a bijection between the
tilings compatible with $(A^*,B^*)$ and the tilings compatible
with some other partition of $\{2,4,...,2k\}$.

Also, recall the variables $x_1,x_2,...,x_m$ that made a brief
appearance in section III before getting swallowed up by the notation.
Leaving aside our appeal to the explicit formulas for 
$s_\sigma(x_1,...,x_k,0,...)$
and $s_\tau(x_1,...,x_k,0,...)$, we may use the linear algebra formalism
of section IV to derive a Schur function identity in infinitely many
variables, expressing the product $s_\sigma s_\tau$ as a sum of products
of other pairs of Schur functions.  It would be desirable to have a
combinatorial explanation of these identities at the level of Young
tableaux.

In section V, we made use of the fact that
if $(A,B)$ is chosen randomly from among the balanced ordered partitions
of $\{1,2,...,2k\}$, and if $s,t \in \{1,2,...,2k\}$ have
the same parity,
then the events $s \in A$ and $t \in A$ are independent of one another.
\We~conjecture, based on numerical evidence, that if $s, t \in \{1,2,...,2k\}$ 
have opposite parity,
then the events $s \in A$ and $t \in A$ are negatively correlated.
This conjecture is made plausible by the fact that the total cardinality of $A$
is required to be $k$.
With the use of this conjecture, one could reduce the bound on the
standard deviation in Theorem 3 by a factor of $\sqrt{2}$.
However, neither Theorem 3 nor this strengthening of it comes
anywhere close to giving a true estimate of the variance of $N_m$,
which empirically is on the order of $\log k$ or perhaps even smaller.

Finally, 
fix $1 \leq k \leq n$.
Define 0-1 random variables $X_1,X_2,...,X_n$
such that for all $(x_1,x_2,...,x_n) \in \{0,1\}^n$,
$\Prob[X_i=x_i \ \mbox{for all $i$}] = 0$
unless $\sum_{i=1}^n x_i = k$,
in which case
\begin{multline*}
\Prob[X_i = x_i \ \mbox{for all $i$}] = \\
\frac{\left( \prod_{1 \leq i < j \leq k} \frac{a_j-a_i}{j-i} \right)
\left( \prod_{1 \leq i < j \leq n-k} \frac{b_j-b_i}{j-i} \right)}{2^{k(n-k)}}\ ,
\end{multline*}
where $\{a_1,a_2,...,a_k\} = \{i:\ x_i = 1\}$  
and $\{b_1,b_2,...,b_{n-k}\} = \{i:\ x_i = 0\}$  
($a_1<a_2<...<a_k$,\ $b_1<b_2<...<b_{n-k}$).
This is the distribution on zig-zag patterns in the $k$th diagonal
of the Aztec diamond, induced by a domino tiling chosen uniformly
at random.
Are the $X_i$'s (nonstrictly) negatively pairwise correlated?

\vspace{0.2in}

\noindent
{\bf Acknowledgment}

We thank Gian-Carlo Rota for his helpful advice and an anonymous
referee for informing us of the information contained in the remark
following the proof of Theorem~2.
Both authors were supported by NSF grant DMS-9206374.
James Propp was supported by NSA grant MDA904-92-H-3060.

\vspace{0.2in}

\pagebreak

\noindent
{\bf Bibliography}

\vspace{0.1in}

\noindent
[\EKLPa]
N.\ Elkies, G.\ Kuperberg, M.\ Larsen, and J.\ Propp,
Alternating sign matrices and domino tilings, part 1,
{\it Journal of Algebraic Combinatorics} {\bf 1},
111-132 (1992).

\vspace{0.1in}

\noindent
[\EKLPb]
N.\ Elkies, G.\ Kuperberg, M.\ Larsen, and J.\ Propp,
Alternating sign matrices and domino tilings, part 2,
{\it Journal of Algebraic Combinatorics} {\bf 1},
219-234 (1992).

\vspace{0.1in}

\noindent
[\fulmek]
M. Fulmek, A Schur function identity, \emph{J.\ Combinatorial Theory
  (A)} \textbf{77}, 177--180 (1997).

\vspace{0.1in}

\noindent
[\Macdo]
I.\ G.\ Macdonald, 
``Symmetric Functions and Hall Polynomials,''
Oxford University Press, Oxford, 1979.

\vspace{0.1in}

\noindent
[\MRR]
W.\ H.\ Mills, D.\ P.\ Robbins, and H.\ Rumsey, Jr.,
Alternating sign matrices and descending plane partitions,
{\it J.\ Combinatorial Theory (A)} {\bf 34},
340-359 (1983).

\vspace{0.1in}

\noindent
[\Sagan]
B.\ E.\ Sagan,
``The Symmetric Group,''
Wadsworth and Brooks/Cole, Pacific Grove, CA, 1991.

\vspace{0.1in}

\noindent
[\Stana]
R.\ Stanley,
Theory and application of plane partitions, Parts 1 and 2,
{\it Studies in Applied Math.} {\bf 50},
167-188, 259-279 (1971).

\vspace{0.1in}

\noindent
[\Stanb]
R.\ Stanley,
``Enumerative Combinatorics,'' volume 2,
Cambridge University Press, New York/Cambridge, 1999.

\vspace{0.1in}

\noindent
[\Turn]
H.\ W.\ Turnbull,
``The Theory of Determinants, Matrices, and Invariants,''
Blackie and Son, London and Glasgow, 1929.

\end{document}